\documentclass[a4paper,11pt]{article}

\usepackage{float}
\usepackage[utf8]{inputenc}
\usepackage[english]{babel}
\usepackage[T1]{fontenc}
\usepackage{amssymb,amsmath,amsthm,mathtools}
\usepackage{mathrsfs}
\usepackage{bm}
\usepackage{color}
\usepackage{url}
\usepackage{hyperref}
\usepackage{comment}
\usepackage{amssymb}
\usepackage{subcaption}

\newcommand{\de}{\partial}

\newcommand{\R}{\mathbb{R}}

\newcommand{\vc}[1]{\boldsymbol{#1}}

\newcommand{\nn}{\mathcal{N}}
\newcommand{\nh}{\hat{\mathcal{N}}}

\newcommand{\red}[1]{\textcolor{red}{#1}}

\theoremstyle{definition}

\usepackage{a4wide}
\usepackage{authblk}

\usepackage[english]{babel}
\usepackage[T1]{fontenc}
\usepackage{amssymb,amsmath,mathtools}
\usepackage{mathrsfs}
\usepackage{bm}
\usepackage{color}
\usepackage{hyperref}
\usepackage{comment}
\usepackage{amssymb}

\usepackage[utf8]{inputenc}
\usepackage{graphicx} 
\usepackage{xcolor}
\usepackage{bbm}
\usepackage{siunitx}

\usepackage{chemformula}
\usepackage{amssymb}
\usepackage{bbold}
\usepackage{minitoc}
\usepackage[normalem]{ulem}
\usepackage{upgreek}

\usepackage{tikz}
\usetikzlibrary{calc}

\newcommand {\dt }{\Delta t}

\newcommand{\new}[1]{\textcolor{black}{#1}}

\title{Variationally mimetic operator network approach to transient viscous flows}

\author[1]{Laura Rinaldi\footnote{laura.rinaldi@imati.cnr.it}}
\author[2]{{Giulio G. Giusteri}\footnote{giulio.giusteri@unipd.it}}
\affil[1]{Consiglio Nazionale delle Ricerche-Istituto di Matematica
Applicata e Tecnologie Informatiche “E. Magenes” (CNR-IMATI), via Ferrata 5/A, 27100 Pavia, Italy}
\affil[2]{Dipartimento di Matematica ``Tullio Levi-Civita'', Universit\`a degli Studi di Padova\\ via Trieste 63, 35131, Padova, Italy}

\date{\today}

\begin{document}

\maketitle

\begin{abstract}
\phantomsection\addcontentsline{toc}{section}{\numberline{}Abstract}

The Variationally Mimetic Operator Network (VarMiON) approach is a machine learning technique, originally developed to predict the solution of elliptic differential problems, that combines operator networks with a structure inherited from the variational formulation of the equations.
We investigate the capabilities of this method in the context of viscous flows, by extending its formulation to vector-valued unknown fields and with a particular emphasis on the space-time approximation context necessary to deal with transient flows.
As a first step, we restrict attention to the regime of low-to-moderate Reynolds numbers, in which the Navier--Stokes equations can be linearized to give the time-dependent Stokes problem for incompressible fluids.
The details of the method as well as its performance are illustrated in three paradigmatic flow geometries where we obtain a very good agreement between the VarMiON predictions and reference finite-element solutions.
\end{abstract}




\section{Introduction}
Machine learning (ML) techniques have been recently applied in many scientific contexts and are showing promising results in forecasting tasks across various fields. 
They are particularly useful either in the presence of incomplete datasets about the physical system, linked to experimental limitations, or when the computational cost of direct numerical solutions of the partial differential equations involved in the mathematical model is excessively high.
Indeed, they allow for a more efficient prediction of new solutions based on a restricted set of training data.

One of the contexts in which we find both experimental limitations and computational challenges is fluid dynamics. For this reason, much effort has been made to develop and apply ML techniques to predict the solutions of flow equations in several situations \cite{brunton2021applying}.
For instance, for predicting solutions of the electroosmotic flow in micro-channels
\cite{baymani2015artificial}, for solving the Reynolds-averaged Navier–Stokes equations for incompressible turbulent flows \cite{eivazi2022physics, ling2016reynolds}, in the aerodynamic flow fields \cite{bhatnagar2019prediction, kou2021data}
as well as in the context of non-Newtonian fluid mechanics \cite{simavilla2025hammering}, or even in cardio-physiology to forecast blood flow \cite{balzotti2022data, garay2024physics, cruz2025enhanced}.

The typical ML techniques \cite{garikipati2024data} developed and used are
Artificial Neural Network (ANN), Convolutional Neural Network (CNN), featuring the two processes of encoding and decoding, 
{Physics-informed} Neural Network (PINN) \cite{raissi2018hidden, raissi2019physics} as well as
Operator Network (ON) \cite{ray2023deep}, especially Deep Operator Network (DeepONet) \cite{ray2023deep, lu2021learning} and Variationally Mimetic Operator Network (VarMiON) \cite{patel2022variationally}.
The latter is an upgrade of DeepONet where the variational formulation of the PDE to predict is exploited.

\begin{figure}[tb]
\centering
\includegraphics[width=0.4\linewidth]{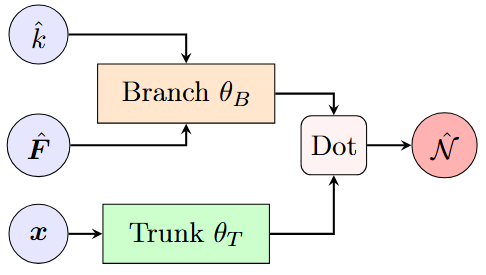}

\caption{Schematic of a generic DeepONet. The output prediction $\hat{\mathcal N}$ is given as a dot product between the output of the trunk (starting from some space input) and the output of the branch network (starting from some input data such as a forcing term). This structure is used to predict the outcome of an operator $\mathcal N: f ,k\mapsto u(x)$.}
\label{fig:deeponet}
\end{figure}

The VarMiON is presented for elliptic equations in the paper \cite{patel2022variationally} (and its further extension in \cite{charles2025optimal}), where the authors describe a new architecture for operator networks that mimics the form of the
numerical solution obtained from an approximation of the variational (or weak) formulation of the
problem. Inspired by the Galerkin approximation method, the
VarMiON is composed of a sub-network, a trunk, whose output constructs the basis functions, and another sub-network, a branch, whose output constructs the coefficients with respect to these basis functions.
The main difference with the conventional DeepONet \cite{ray2023deep}, is that the architecture of the branch sub-network is precisely determined by the discrete weak form \cite{patel2022variationally}. 
Especially, let us consider the elliptic PDE
\begin{align*}
        -k\Delta  u&= f(\vc x) \qquad \vc x \in \Omega,\\
        u(\vc x)&=g(\vc x)  \qquad \vc x \in \partial\Omega,
\end{align*}
and define the operator $\mathcal{N}$ that maps the data $k$ and $f$ to the solution $u$ of the PDE
$$u = \mathcal N(f,k)(\vc x)$$
for $g$ given.
In this example, the DeepONet comprises two sub-networks that construct the basis functions and the coefficients for the output to approximate ${\mathcal{N}}$ (Figure \ref{fig:deeponet}). The approximated operator $\hat {\mathcal{N}}$ maps the sampled data $\hat k,$ $\hat{\vc F}$ to the prediction $\hat u$ computed by the network: $$\hat{u}=  \hat {\mathcal{N}}(\hat{\vc F},\hat{k}).$$

The VarMiON approximates $ {\mathcal{N}}$ by means of two sub-networks as well, but the architecture of the branch is determined by the weak form of the PDE, as detailed in Section \ref{sec:varmion} below. 
To our knowledge, the variational formulation is already used in ML to inform the network, as done for VarNet and VPINNs \cite{khodayi2020varnet, kharazmi2019variational}, but only to modify the loss function and not the structure of the net, as is done in the VarMiON scheme.

In this paper, we focus our attention on fluid dynamics, particularly the time-dependent Stokes equation, as a first step to address the more challenging Navier--Stokes system.
Much effort has been made from the point of view of ML (as reviewed in \cite{brunton2021applying, lakshmiexploring, kutz2017deep, taira2025machine}) to predict the outcome of the velocity and the pressure field.
As mentioned above, ANN represents a possible way to attack the problem of predicting the solution of the Navier--Stokes equation \cite{baymani2015artificial, mccracken2018artificial}.  
Further methods involve PINNs  \cite{eivazi2022physics,raissi2018hidden, amalinadhi2022physics} or CNN  \cite{miyanawala2017efficient, ma2021physics,ranade2021discretizationnet}.

Our idea is to use the VarMiON technology for time-dependent equations, as recently done for scalar elliptic problems \cite{rinaldi2025, Rinaldi_Marcuzzi_Chinellato_dts}, and to adapt it for predicting a vector function that depends on space-time variables. 
We restrict attention to the regime of low-to-moderate Reynolds numbers, in which the Navier--Stokes equations can be linearized to give the time-dependent Stokes problem for incompressible fluids.
In Section \ref{sec:varmion} we describe the adaptation of the VarMiON method to our context.
In Section \ref{sec:num_results}, we show the numerical results obtained by training the VarMiON for the paradigmatic examples of the cavity flow, the flow past a cylinder, and the contraction flow, with a particular emphasis on the transient behavior.
We summarize our results and mention further research directions in Section~\ref{sec:conclusions}.

\section{VarMiON formulation}\label{sec:varmion}

In this section, we illustrate how the Variationally Mimetic Operator Networks (VarMiON) presented in the paper \cite{patel2022variationally} (and its further extension in \cite{charles2025optimal}), can be used to predict the solution of fluid-dynamics problems.
VarMiON is an operator network (ON) that tries to forecast the solution of a PDE given as the output of an operator $\mathcal{N}$, when the input data of $\mathcal{N}$ are initial conditions, boundary or forcing terms, and parameters of the differential problem. 

Like the conventional DeepONet \cite{ray2023deep}, the
VarMiON is composed of two types of sub-networks, the \emph{trunks}, that construct the basis functions for the output, 
and the \emph{branches}, that construct the coefficients for the unknown reconstruction. The peculiarity of the method is that the architecture of the branches is determined by the discrete weak form of the PDE \cite{patel2022variationally}. 

\new{Both VarMiON (and DeepONet) and classical reduced-order (ROM) techniques (such as projection- or singular-value-based approaches) aim to represent the solution vector as a linear combination of coefficients and basis functions. 
In standard reduced-order modeling, these basis functions—often referred to as ``modes'', and the associated ``modal coefficients'', are derived by projecting the original PDE system onto a lower-dimensional, computationally tractable subspace.}

\new{In particular, the modes are low-dimensional spatial functions constructed from a matrix $\vc{\mathcal U}$ whose columns are solution snapshots. From this snapshot matrix, modes are extracted via methods like singular value decomposition (SVD) or by factorizing the correlation matrix. The modal coefficients are then obtained by projecting the solution snapshots onto these modes.}

\new{In contrast, in operator-learning frameworks such as DeepONet and VarMiON, the branch network is trained to encode the input function into coefficients of an implicit (latent) basis, while the trunk network is trained to produce the corresponding basis functions over the space-time domain. These two components are subsequently combined to yield the solution.} 

\new{Concretely, the branch network takes the input function and maps it to a latent feature representation. The trunk network takes the evaluation point $(\vc x, t)$ and produces another latent feature representation. The operator’s value at $(\vc x, t)$ is then computed as the inner product between these two latent representations. Both the trunk and branch components are realized as neural networks.  This approach learns an optimal set of basis functions directly from data, enabling it to capture complex solution manifolds. Furthermore, the data-driven basis and coefficients remove the need for hand-crafted basis functions, which are typically required in traditional numerical schemes \cite{choi2026operator}.}
In the case of VarMiON, the global vectors of trainable parameters are the combination of the ones of the branches and of the trunks.

\subsection{The flow equations}
The planar Navier--Stokes problem that describes the behavior of Newtonian viscous fluids in a domain $\Omega\subseteq\R^2$ with Lipschitz boundary and on a time interval $[0,\tau]$ consists of a system of equations for the velocity $\vc u= (u_1, u_2)^\top $ and pressure $p$. In the incompressible case, with constant and uniform mass density $\rho>0$, it reads
\begin{equation}\label{eq_Navier--Stokes}
\begin{aligned}
    \rho\left( \frac{\partial \vc u}{\partial t} + (\vc u \cdot \nabla) \vc u\right) = \nabla \cdot \vc\sigma(\vc u, p) + \vc f,\\
    \nabla \cdot \vc u =0,
\end{aligned}
\end{equation}
where $\vc f=(f_1, f_2)^\top$ is the body force per unit volume and $\vc\sigma(\vc u, p)$ denotes the stress tensor which, for a Newtonian fluid, is given by
\begin{equation}
    \vc\sigma(\vc u, p) = 2 \mu \dot{\vc\varepsilon}(\vc u ) -p \vc I,
\end{equation}
with $\vc I$ the identity tensor, $\mu>0$ the dynamic viscosity, and $\dot{\vc\varepsilon}(\vc u)$ the strain-rate tensor defined as
\begin{equation}
    \dot{\vc\varepsilon}(\vc u) := \frac{1}{2}\left ( \nabla\vc u + (\nabla \vc u)^{\top} \right).
\end{equation}

When inertial forces, that increase with the magnitude of $\vc u$, are small compared to viscous forces, the Navier--Stokes equation \eqref{eq_Navier--Stokes} can be linearized to give, in the time interval $[0,\tau]$, the time-dependent Stokes problem on $\Omega\times [0,\tau]$ as
\begin{equation}\label{eq:stokes}
\begin{aligned}
    \rho \frac{\partial \vc u}{\partial t}  = -\nabla p+\mu\Delta \vc u + \vc f,\\
    \nabla\cdot\vc u=0,
\end{aligned}
\end{equation}
in which we substituted the Newtonian form of the stress tensor.
The differential system~\eqref{eq:stokes} must be accompanied by initial conditions on $\Omega \times \{0\}$, given by
\begin{equation}
    \vc u(\vc x, 0) = \vc u_0(\vc x),\qquad
    p(\vc x, 0) = p_0(\vc x), 
\end{equation}
and boundary conditions, that will be specified for each of the problems solved in Section~\ref{sec:num_results} and for now we summarize in a vector field $\vc g$ defined on $\de\Omega\times[0,\tau]$.

\subsection{Solution operators and input data}

We introduce the space $X$ of the input data represented by $(\rho, \mu, \vc f, \vc u_0, p_0, \vc g )$ by defining
\[
X:=\R  \times \R  \times L^2([0, \tau], L^2(\Omega,\R^2)) \times L^2(\Omega,\R^2) \times L^2(\Omega) \times L^2([0, \tau], H^\frac{3}{2}(\de\Omega,\R^2))
\]
and the space $V$ of solutions $(\vc u,p)$ as
\[
V:=L^2([0,\tau],H^1_0(\Omega,\R^2))\times L^2([0,\tau],L^2(\Omega)).
\]
We can thus consider the solution operator $\mathcal{N}:X \to V$ that maps the input data to the solution of the Stokes problem.
The goal of the VarMiON is to learn a good approximation $\nh$ of the solution operator $\nn$.

The input data for the VarMiON solution operator $\nh$ are obtained by applying a so-called \emph{sensing operator} $\new{\hat{\mathcal{Q}}}$ that evaluates the input functions of $\mathcal{N}$ at some time nodes $\{ \hat{t}_i\}_{i=1}^r$, some space nodes $\{ \hat{\vc x}_i\}_{i=1}^{k}$ and, some boundary nodes $\{ \hat{\vc x}_i\}_{i=1}^{k'}$, mimicking what could be measured by idealized sensors. We have
\begin{equation}
\new{\hat{\mathcal{Q}}}:\left\{
\begin{aligned}
X\qquad\qquad&\to \R\times\R\times\R^{2kr}\times\R^{2k}\times\R^{k}\times\R^{2k'r}\\
( \rho, \mu,\vc f,\vc u_0, p_0,\vc g) &\mapsto \qquad (\rho,\mu, \hat{\vc {F}} , \hat{\vc U}_0, \hat{\vc P}_0, \hat{\vc G})
\end{aligned}
\right.
\end{equation}
where the components of the vectors $\hat{\vc{F}}$, $\hat{\vc{U}}_0$, $\hat{\vc{P}}_0$, and $\hat{\vc{G}}$ contain the nodal values of the corresponding fields.
The sensed data $\new{\hat{\mathcal{Q}}}(\rho, \mu,\vc f,\vc u_0, p_0,\vc g)$ will be given as inputs to the approximate solution operator $\nh$.

The general regularity of the data may not force point-wise values of the functions, but we can use Lebesgue values or regularized data to make our procedure uniquely defined.

\subsection{Space-time variational formulation}
To take advantage of the structure of VarMiON presented for the elliptic equations \cite{patel2022variationally}, we work with the weak space-time formulation that discretizes the space and time variables simultaneously \cite{rinaldi2025, Rinaldi_Marcuzzi_Chinellato_dts,frank2021implementation, langer2021space, von2023time}.

By introducing test fields $\vc v$ and $q$ on $\Omega \times [0,\tau]$ for velocity and pressure, respectively, the weak form of the Stokes problem~\eqref{eq:stokes}, easily obtained by standard computations, reads
\begin{multline}\label{eq:weak-con-pb}
\int_0^\tau\int_{\Omega}\left(\rho\frac{\de\vc u}{\de t}\cdot\vc v+q\nabla\cdot\vc u-p\nabla\cdot\vc v+\mu\nabla\vc u\cdot\nabla\vc v-\vc f\cdot\vc v\right)\,d\vc x dt\\
+\int_0^\tau\int_{\de\Omega}\left(p\vc n-\mu\frac{\de\vc u}{\de n}\right)\cdot\vc v\,dS dt=0
\end{multline}
for any pair $(\vc v,q)$ of test fields in suitably defined spaces adapted to each specific example.

We employ the Galerkin method by defining two types of basis functions for $t>0$:
\begin{itemize}
    \item $\pi_i(\vc x, t)$ corresponding to the pressure;
    \item $\vc\psi_{ij}(\vc x, t)$ corresponding to the $j$-th velocity component, 
where the velocity basis functions in vector form are
\begin{equation}
    \vc \psi_{i1} = \begin{bmatrix}
        \psi_i\\
        0
    \end{bmatrix} 
    \quad\text{and}\quad
    \vc \psi_{i2} = \begin{bmatrix}
        0\\
        \psi_i
    \end{bmatrix}.
\end{equation}
\end{itemize}

Now, the approximate functions can be written as
\begin{align}
    &\vc u_h = \sum_j^{N_h} u_{1j} \vc\psi_{j1} + u_{2j} \vc\psi_{j2} = \sum_j^{2N_h} u_j \vc \psi_j = \vc{U}^\top \vc{\Psi},\\
        &p_h = \sum_{j}^{M_h} p_j \pi_j = \vc{P}^\top \vc{\Pi}.
\end{align}
We group in the vector $\vc{\Phi}=(\vc\Psi,\vc\Pi)^\top$ all of the basis functions.

We denote by $\new{{\mathcal{Q}}}$ the projector on the finite-dimensional subspace associated with the Galerkin approximation, namely $\new{{\mathcal{Q}}}( \rho, \mu,\vc f, \vc u_0, p_0, \vc g)=( \rho, \mu,\vc f_h, \vc u_{0_h}, p_{0_h}, \vc g_h)$, and we obtain the discrete matrix formulation of problem \eqref{eq:weak-con-pb} in the form
\begin{equation}\label{eq:matrix_da}
\begin{aligned}
       \mathbf{L}\vc{U}=\vc 0,\\
       (\mathbf{W}(\rho)\vc{U}+   \mathbf{K}(\mu)\vc{U}) -\mathbf{L}^{\top}\vc{P} =  \mathbf{M} \vc{F},
\end{aligned}
\end{equation}
where
\begin{align}
       \mathbf{L} &= (\pi_i, \nabla\cdot \vc\psi_j),\quad
       \mathbf{W}(\rho)& =(\rho \partial_t \vc\psi_i , \vc\psi_j),\quad
       \mathbf{K}(\mu) = (2\mu \nabla\vc\psi_i, \nabla\vc\psi_j),\quad
       \mathbf{M} = (\vc\psi_i, \vc\psi_j).
\end{align}

We set $\mathbf A(\rho, \mu) :=\mathbf{W}(\rho)  + \mathbf{K}(\mu)$ and consider the decompositions
\[
\vc U=(\vc U_0,\vc G,\tilde{\vc U})^\top,
\quad
\vc P=(\vc P_0,\tilde{\vc P})^\top,
\quad
\mathbf{A}=(\mathbf{A}_0|\mathbf{A}_g|\tilde{\mathbf{A}}),
\quad
\text{and}
\quad
\mathbf{L}=(\mathbf{L}_0|\tilde{\mathbf{L}}),
\]
where we separated the degrees of freedom pertaining to initial conditions (subscript $0$) and boundary conditions ($\vc G$) from the remaining components (denoted with a superimposed tilde). 
With these provisions, we rewrite the problem \eqref{eq:matrix_da} as

\begin{equation}
    \begin{bmatrix}
    \tilde{\mathbf{A}} & -\tilde{\mathbf{L}}^{\top}\\
    \tilde{\mathbf{L}} & \vc 0
\end{bmatrix}
\begin{bmatrix}
   \tilde{\vc{U}}\\
    \tilde{\vc P}
\end{bmatrix}
=
\begin{bmatrix}
    \mathbf{M} \vc{F} -\mathbf{A}_{g} {\vc{G}} - \mathbf{A}_0 \vc{U}_0 + \mathbf{L}^{\top}_0\vc{P}_0\\
    -  \mathbf{L}_0\vc{U}_0
\end{bmatrix}.
\end{equation}
Notice that, the basis $ \vc \psi_j(\vc x, t)$ and $\pi_j(\vc x, t)$ are chosen such that the matrix on the left-hand side is invertible for any $\rho$ and $\mu$.
From the previous equation we then obtain
\begin{align}
\begin{bmatrix}
   \tilde{\vc{U}}\\
    \tilde{\vc P}
\end{bmatrix}
&=    \begin{bmatrix}
    \mathbf{D}_1 & \mathbf{D}_2\\
    \mathbf{D}_3 & \mathbf{D}_4
\end{bmatrix}
\begin{bmatrix}
    \mathbf{M} \vc{F} -\mathbf{A}_{g} {\vc{G}} - \mathbf{A}_0 \vc{U}_0 + \mathbf{L}^{\top}_0\vc{P}_0\\
    -  \mathbf{L}_0\vc{U}_0
\end{bmatrix} \nonumber \\
&=    \begin{bmatrix}
    \mathbf{D}_1 & \mathbf{D}_2\\
    \mathbf{D}_3 & \mathbf{D}_4
\end{bmatrix}
\left(
\begin{bmatrix}
    \mathbf{M} \vc{F} \\
    0
\end{bmatrix}
+
\begin{bmatrix}
    -\mathbf{A}_{g} {\vc{G}} \\
    0
\end{bmatrix}
+
\begin{bmatrix}
 - \mathbf{A}_0 \vc{U}_0 \\
    -  \mathbf{L}_0\vc{U}_0
\end{bmatrix}
+
\begin{bmatrix}
     \mathbf{L}^{\top}_0\vc{P}_0\\
    0
\end{bmatrix}
\right),\label{eq:fem_matrix}
\end{align}
where $\mathbf{D}_i$, $i=1,\ldots,4$, denote suitable blocks of the inverse matrix.
The discrete solution operator $ {\mathcal{N}_h}$ for $t>0$ is defined as
\begin{multline}\label{eq:discrete_op}
    {\mathcal{N}_h}(\rho, \mu, \vc f_h,  \vc {u}_{0_h} ,  p_{0_h}, \vc g_h)= 
    (\vc u_h(\cdot, \cdot, \rho, \mu, \vc f_h , \vc{ u}_{0_h} ,  p_{0_h}, \vc g_h), p_h(\cdot, \cdot, \rho, \mu, \vc f_h , \vc{ u}_{0_h} ,  p_{0_h}, \vc g_h))^\top\\
    =  \left(    \mathbf{b}_f(\rho, \mu, \vc f_h) +  \mathbf{b}_g (\rho, \mu, \vc g_h)+  \mathbf{b}_{u0} (\rho, \mu, \vc u_{0_h}) +  \mathbf{b}_{p0} (\rho, \mu, p_{0_h})\right)^{\top} 
       { \vc \Phi} ,
\end{multline}
where 
\begin{align*}
    &\mathbf{b}_f(\rho, \mu, \vc f_h)=  \begin{bmatrix}
    \mathbf{D}_1 \mathbf{M} \vc{F}\\
    \mathbf{D}_3 \mathbf{M} \vc{F} 
\end{bmatrix} ,\\
    & \mathbf{b}_g (\rho, \mu, \vc g_h)= - \begin{bmatrix}
    \mathbf{D}_1 \mathbf{A}_{g} {\vc{G}} \\
    \mathbf{D}_3 \mathbf{A}_{g} {\vc{G}} 
\end{bmatrix} ,\\
& \mathbf{b}_{u0} (\rho, \mu, \vc{ u}_{0_h})=  - \begin{bmatrix}
    \mathbf{D}_1 \mathbf{A}_0 \vc{U}_0 \\
    \mathbf{D}_3 \mathbf{A}_0 \vc{U}_0 
\end{bmatrix} - \begin{bmatrix}
     \mathbf{D}_2 \mathbf{L}_0 \vc{U}_0 \\
     \mathbf{D}_4 \mathbf{L}_0 \vc{U}_0
\end{bmatrix}, \\
    & \mathbf{b}_{p0} (\rho, \mu ,  p_{0_h})=  \begin{bmatrix}
    \mathbf{D}_1\mathbf{L}^{\top}_0\vc{P}_0\\
    \mathbf{D}_3\mathbf{L}^{\top}_0\vc{P}_0
\end{bmatrix} .
\end{align*}

\subsection{VarMiON structure}

The VarMiON mimics the structure of ${\mathcal{N}_h}$ starting from the input data that we presented at the beginning of this section. The final VarMiON operator $\hat{\mathcal{N}}$ is given as a dot product between the output of the branches
\begin{align}
    \begin{bmatrix}
    \hat{\vc U}\\
    \hat{P}
\end{bmatrix} ={\mathbf D}({\rho}, {\mu})(\new{\hat{\mathbf A}_{F}}\hat{\vc F} -\new{\hat{\mathbf A}_{G}}\hat{\vc G}-\new{\hat{\mathbf A}_{U_0}}\hat{\vc U}_0+ \new{\hat{\mathbf A}_{P_0}}\hat{\vc P}_0)  \label{eq:varmion_matrix}
\end{align}
and the output of three trunks $\vc{\uptau}$, one for each component of the solution field.
We have
\begin{multline}\label{eq:varmion_op}
    \hat{\mathcal{N}} (\rho, \mu,\hat{\vc F},\hat{\vc U}_0,\hat{\vc P}_0, \hat{\vc G})=(\hat{\vc u}(\cdot, \cdot,  \rho, \mu,\hat{\vc F},\hat{\vc U}_0,\hat{\vc P}_0, \hat{\vc G}), \hat{p}(\cdot, \cdot,  \rho, \mu,\hat{\vc F},\hat{\vc U}_0,\hat{\vc P}_0, \hat{\vc G}) )^\top
    \\
    =  \left(    \vc{\upbeta}_f(\rho, \mu,\hat{\vc F})  +\vc{\upbeta}_g (\rho, \mu,\hat{\vc G})+ \vc{\upbeta}_{u0} (\rho, \mu,\hat{\vc U}_0) + \vc{\upbeta}_{p0} (\rho, \mu,\hat{\vc P}_0)\right)^{\top}\vc \uptau .
\end{multline}
Notice that \eqref{eq:varmion_op} replicates the structure of \eqref{eq:discrete_op}.

\begin{figure}[H]
\includegraphics[width=\linewidth]{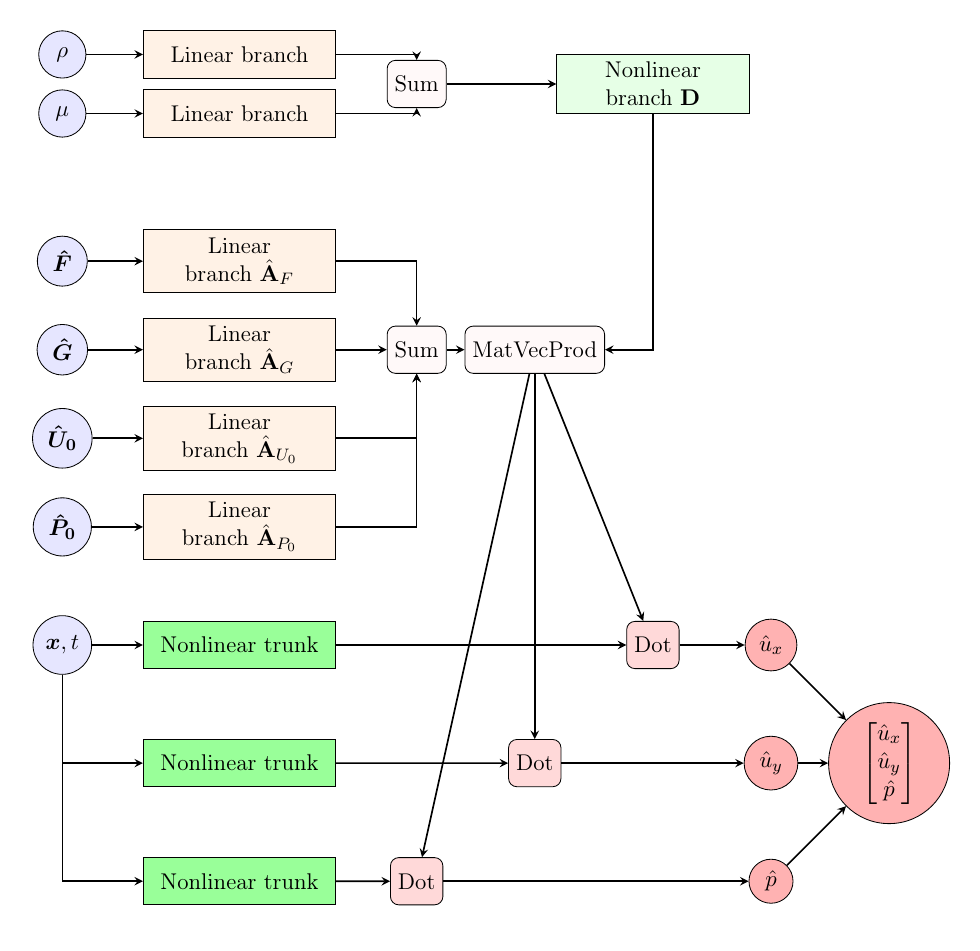}

\caption{VarMiON scheme for the time-dependent Stokes equation. The network prediction $\hat{\mathcal{N}}$ is given as a dot product between the outcome of the trunks and the outcome of the branches starting from the input data $(\rho, \mu,\hat{\vc F},\hat{\vc U}_0,\hat{\vc P}_0, \hat{\vc G})$.}
\label{fig:VarMiON_ns}
\end{figure}
We consider the VarMiON shown in Figure \ref{fig:VarMiON_ns} and summarized in Table \ref{tab:fem_var}, which comprises:
\begin{itemize}
    \item  Two linear branches taking the input $ {\mu}$ and $ {\rho}$, which sum is transformed into a matrix output ${\mathbf D}(\rho, \mu) \in \mathbb R^{p\times p}$. Here $p$ is the latent dimension of the VarMiON.
    \item A linear branch taking the input $\hat{\vc F}$ and transforming it as $\new{\hat{\mathbf A}_{F}}\hat{\vc F} $, where $\new{\hat{\mathbf A}_{F}} \in \mathbb R^{p\times 2kr}$ is a learnable matrix. The output of this branch is acted upon by the matrix ${\mathbf D}$ to give ${\vc{\upbeta}}_f(\rho,\mu,\hat{\vc F} )={\mathbf D}(\rho, \mu)\new{\hat{\mathbf A}_{F}}\hat{\vc F} \in \mathbb R^p$.
    \item A linear branch taking the input $\hat{\vc G}$ and transforming it as ${\new{\hat{\mathbf A}_{G}}}\hat{\vc G} $, where ${\new{\hat{\mathbf A}_{G}}} \in \mathbb R^{p\times 2k'r}$ is a learnable matrix. The output of this branch is acted upon by the matrix ${\mathbf D}$ to give ${\vc{\upbeta}_g}( \rho, \mu,\hat{\vc F} )={\mathbf D}(\rho, \mu)\new{\hat{\mathbf A}_{G}}\hat{\vc G} \in \mathbb R^p$.
    \item A linear branch taking the input $\hat{\vc U}_0$ and transforming it as ${\new{\hat{\mathbf A}_{U_0}}}\hat{\vc U}_0 $, where $\new{\hat{\mathbf A}_{U_0}} \in \mathbb R^{p\times 2k}$ is a learnable matrix. The output of this branch is acted upon by the matrix ${\mathbf D}$ to give $\vc{\upbeta}_{u0}( \rho,\mu,\hat{\vc U}_0 )={\mathbf D}(\rho, \mu)\new{\hat{\mathbf A}_{U_0}}\hat{\vc U}_0 \in \mathbb R^p$.
    \item A linear branch taking the input $\hat{\vc P}_0$ and transforming it as ${\new{\hat{\mathbf A}_{P_0}}}\hat{\vc P}_0 $, where $\new{\hat{\mathbf A}_{P_0}} \in \mathbb R^{p\times k}$ is a learnable matrix. The output of this branch is acted upon by the matrix ${\mathbf D}$ to give $\vc{\upbeta}_{p0}( \rho, \mu,\hat{\vc P}_0 )={\mathbf D}(\rho, \mu)\new{\hat{\mathbf A}_{P_0}}\hat{\vc P}_0 \in \mathbb R^p$.
    \item Three non-linear trunks taking the input $(\vc x, t) \in \Omega \times (0, \tau]$, which give the output 
    
    $\vc \uptau(\vc x, t) =
(\uptau_1(\vc x, t), \dots , \uptau_p(\vc x, t))^{\top}$, where each $\uptau_i : \Omega \times(0, \tau] \rightarrow \mathbb R$ is a trainable network.
\end{itemize}
\begin{table}
\centering
\begin{tabular}{|p{4cm}|p{4cm}|p{3cm}|p{3cm}|}
\hline
\textit{FEM input}& \textit{FEM terms} & \textit{VarMiON terms } & \textit{Sub-networks}\\
\hline 
\begin{center}
${\mu}$, ${\rho}$, pde structure 
\end{center}& $$\begin{bmatrix}
    \mathbf{D}_1 & \mathbf{D}_2\\
    \mathbf{D}_3 & \mathbf{D}_4
\end{bmatrix}$$ & $$\red{\mathbf D}(\hat{\vc {\rho}},\hat{\vc {\mu}})$$& \begin{center}$2$ linear branches\end{center} \\
\hline
\begin{center} external force $ \vc{F}$ \end{center}& $$\begin{bmatrix}
    \mathbf{D}_1 \mathbf{M} \vc{F}\\
    \mathbf{D}_3 \mathbf{M} \vc{F} 
\end{bmatrix}$$ & $${\mathbf D}({\rho}, {\mu})\red{{\hat{\mathbf A}_{F}}}\hat{\vc F} $$& \begin{center}
    linear branch ${\hat{\mathbf A}_{F}}$
\end{center} \\
\hline
\begin{center}boundary data $ \vc{G}$ \end{center}& $$ \begin{bmatrix}
    \mathbf{D}_1 \mathbf{A}_{g} {\vc{G}} \\
    \mathbf{D}_3 \mathbf{A}_{g} {\vc{G}} 
\end{bmatrix} $$& $${\mathbf D}({\rho}, {\mu})\red{{\hat{\mathbf A}_{G}}}\hat{\vc G}$$  & \begin{center}linear branch ${\hat{\mathbf A}_{G}}$\end{center} \\
\hline
\begin{center}initial velocity $ \vc{U}_0$ \end{center}& $$ \begin{bmatrix}
    \mathbf{D}_1 \mathbf{A}_0 \vc{U}_0 \\
    \mathbf{D}_3 \mathbf{A}_0 \vc{U}_0 
\end{bmatrix} + \begin{bmatrix}
     \mathbf{D}_2 \mathbf{L}_0 \vc{U}_0 \\
     \mathbf{D}_4 \mathbf{L}_0 \vc{U}_0
\end{bmatrix}$$  & $${\mathbf D}({\rho}, {\mu})\red{{\hat{\mathbf A}_{U_0}}}\hat{\vc U}_0$$ & \begin{center}linear branch ${\hat{\mathbf A}_{U_0}}$\end{center} \\
\hline
\begin{center}initial pressure $ \vc{P}_0$\end{center}& $$
    \begin{bmatrix}
    \mathbf{D}_1\mathbf{L}^{\top}_0\vc{P}_0\\
    \mathbf{D}_3\mathbf{L}^{\top}_0\vc{P}_0
\end{bmatrix}$$ & $${\mathbf D}({\rho}, {\mu})\red{{\hat{\mathbf A}_{P_0}}}\hat{\vc P}_0$$ & \begin{center}linear branch ${\hat{\mathbf A}_{P_0}}$\end{center} \\
\hline
\begin{center}space-time variables $ x, t$\end{center}& $$\vc\Phi$$ & $$\red{\vc\uptau}$$ & \begin{center}non-linear trunks \end{center}\\
\hline
\end{tabular}
\caption{FEM and VarMiON correspondence. Each sub-network learns the matrix in red.}
\label{tab:fem_var}
\end{table}

\subsection{VarMiON samples and loss}
The discrete operator ${\mathcal{N}}_h$ is used to generate samples for training VarMiON in the following manner:
\begin{enumerate}
    \item  For $1 \leq j \leq J$, consider distinct samples of the data $(\rho^j, \mu^j, \vc f^j, \vc u_0^j, p_0^j, \vc g^j) \in X$.
    \item Use the projector $\new{{\mathcal{Q}}}$ to obtain the discrete approximations  $(\rho^j, \mu^j, \vc f^j_h, \vc u_{0_h}^j, p_{0_h}^j, \vc g_h^j)$.
    \item Find the discrete numerical solution $(\vc u^{j}_h, p_h^j)^\top= \mathcal{N}_h(\rho^j, \mu^j, \vc f^j_h, \vc u_{0_h}^j, p_{0_h}^j, \vc g_h^j).$
    \item Use the $\new{\hat{\mathcal{Q}}}$ to generate the VarMiON input vectors $( \rho^j, \mu^j,\hat{\vc F}^j,\hat{\vc U}_0^j,\hat{\vc P}_0^j, \hat{\vc G}^j)$.
    \item Select a set of output nodes $\{ \vc {x}_l, t_i\}_{l=1}^L,_{i=1}^M$. For each $1 \leq j \leq J$, sample the numerical solution at these nodes as $(\vc u^{j}_h(\vc {x}_l, t_i), p^{j}_h(\vc {x}_l, t_i))^\top$.
    \item Compose the training set containing $J\times L \times M$ samples 
    \begin{equation}
    S=\{\underbrace{{ \rho^j,\mu^j,\hat{\vc F}^j,\hat{\vc U}_0^j,\hat{\vc P}_0^j, \hat{\vc G}^j}}_{\text{branch inputs}} \quad , \underbrace{{ \vc {x}_l, t_i}}_{\text{trunk inputs}},  \ \underbrace{{(\vc u^{j}_h(\vc {x}_l, t_i), p^{j}_h(\vc {x}_l, t_i))^\top}}_{\text{target output}}\}.
    \end{equation}
\end{enumerate}

To train the VarMiON, we have to optimize a loss function $\Pi$ computed as the $L_2$ norm of the difference between the solution and the prediction computed on the output nodes $\{\vc{x}_l, t_i\}$
\begin{equation}
  \sum_{j=1}^{J}
\int_{0}^{\tau}\int_{\Omega} ((\vc u^{j}_h(\vc {x}, t), p^{j}_h(\vc {x}, t))^\top- \hat{\mathcal{N}}(\rho^j, \mu^j,\hat{\vc F}^j,\hat{\vc U}_0^j,\hat{\vc P}_0^j, \hat{\vc G}^j)( \vc x, t) )^2 \, d\vc x \, dt .
\end{equation}

We choose $\{ \omega_{li}:  1 \le l \le L , 1 \le i \le M \} \subset \mathbb{R}$ (positive) weights, corresponding to the nodes $\{(\vc {x}_l, t_i)\}_{l=1}^L,_{i=1}^M$, for the numerical integration of the square of a function on $(0, \tau] \times \Omega$, i.e. such that for a function $\xi: (0, \tau] \times \Omega \to \mathbb{R}$:
\begin{equation}
   \left |\sum_{l=1}^L\sum_{i=1}^M \omega_{li} \xi^2( \vc {x}_l, t_i) -\int_0^{\tau} \int_{\Omega}\xi^2( \vc {x},t) dt d\vc {x} \right| \leq \frac{C(\xi)}{L^{\gamma}}
\end{equation}
where the constant $C$ depends on $\xi$ and, $\gamma$ is the rate of convergence that depends on the quadrature rule used.

Denoting by $\vc\theta$ the vector of trainable parameters of our VarMiON, we have
\begin{align}\label{error_l2}
    &\Pi(\vc \theta)=\frac{1}{J}\sum_{j=1}^J \Pi_j(\vc \theta),\\
    &\Pi_j(\vc \theta)= \sum_{l=1}^L\sum_{i=1}^M \omega_{li}((\vc u^{j}_h(\vc {x}_l, t_i), p^{j}_h(\vc {x}_l, t_i))^\top-\hat{\mathcal{N}}_{\vc \theta}(\rho^j,\mu^j,\hat{\vc F}^j,\hat{\vc U}_0^j,\hat{\vc P}_0^j ,\hat{\vc G}^j)(\vc {x}_l, t_i) )^2.
\end{align}
Notice that
\begin{equation}
    \Pi_j(\vc \theta) \approx \| \mathcal{N}_h \circ  {\new{{\mathcal{Q}}}}(\rho^j, \mu^j, \vc f^j, \vc u_0^j, p_0^j, \vc g^j) -  \hat{\mathcal{N}}_{\vc \theta} \circ   \new{\hat{\mathcal{Q}}}(\rho^j, \mu^j, \vc f^j, \vc u_0^j, p_0^j, \vc g^j)\|^2_{L_2}
\end{equation}
and, training the VarMiON corresponds to solving the optimization problem 
\begin{equation}
    \vc \theta^*=\underset{\vc \theta}{\mathrm{argmin}} \, \Pi(\vc \theta).
\end{equation}

\section{Numerical results}\label{sec:num_results}
In this section, we show the performance of the VarMiON for three different case studies: the cavity flow, the flow past a cylinder, and the contraction flow. 
To this end, we first compute FEM solutions by randomly sampling the values of the viscosity $\mu$ and of the external forcing terms $\vc f$ in a prescribed interval. A portion of the data is then used to train the VarMiON, while the remaining data are used to test the quality of the VarMiON prediction.
\begin{figure}[t]
    \centering
    \begin{minipage}{0.33\linewidth}
    (a)
    \includegraphics[width=1.0\linewidth]{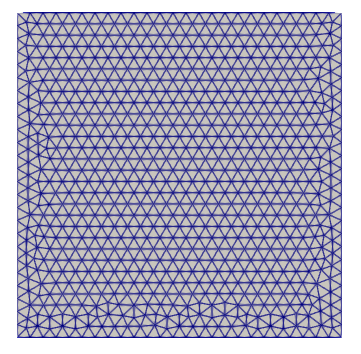}
    \end{minipage}%
    \begin{minipage}{0.675\linewidth}
    (b)
    
    \includegraphics[width=\linewidth]{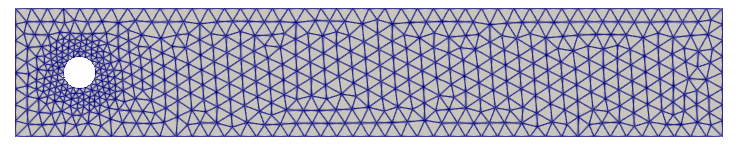}
    
    (c)\vspace{-2mm}
    
    \includegraphics[width=\linewidth]{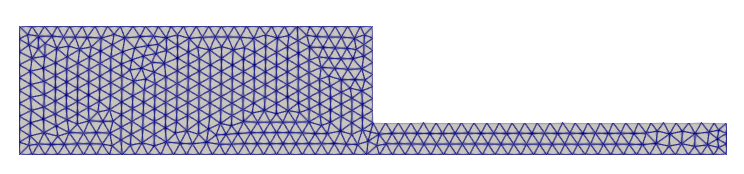}
    \end{minipage}
    
    \caption{The meshes used for the FEM solution of (a) the cavity flow, (b) the flow past a cylinder, and (c) the contraction flow.}
    \label{fig:mesh}
\end{figure}
The FEM solutions are obtained with a modified version of Chorin’s method, the so-called incremental pressure correction scheme, drawing on an implementation in the FEniCSx environment \cite{online_ns_chorin}.

The finite-element meshes used for each example are shown in Figure \ref{fig:mesh}.
Each FEM solution is sampled on a set of points distributed on a rectangular lattice within the space domain repeated on a discrete set of temporal points. The VarMiON is trained to predict the values of the unknown fields on that space-time lattice. \new{In Table \ref{tab_details} we summarize the details of the learning procedure for each test case.}

The structure of the VarMiON network used for these test cases is a simplified scheme of the one in Figure \ref{fig:VarMiON_ns}. We keep only the branches starting from $\mu$ and $\vc f$ \new{because the density $\rho$ is fixed, equal to one, and initial conditions are assumed to be zero, so they are not learned.}
\new{The simplified architecture leads to a reduced formulation of equation \eqref{eq:varmion_op} as given by
\begin{multline}
    \hat{\mathcal{N}} (\mu,\hat{\vc F},\hat{\vc U}_0,\hat{\vc P}_0, \hat{\vc G})=(\hat{\vc u}(\cdot, \cdot, \mu,\hat{\vc F}, \hat{\vc G}), \hat{p}(\cdot, \cdot, \mu,\hat{\vc F}, \hat{\vc G}) )^\top= \left(    \vc{\upbeta}_f( \mu,\hat{\vc F})  +\vc{\upbeta}_g ( \mu,\hat{\vc G})\right)^{\top}\vc \uptau .
\end{multline}}
\new{In particular, as will be discussed in the following subsections, for the cavity flow test the forcing term is set to zero and the boundary condition is a velocity supplied as input to the network. Conversely, in the flow past a cylinder and in the contraction flow, the homogeneous boundary conditions are neglected by the network architecture and the forcing term, generated by a pressure gradient, is taken into account.}
\new{In the preceding section, we presented a more general formulation that is ready to be applied in a variety of situations.}

The operator architecture of each branch follows the one originally presented for elliptic PDEs (see the Appendix in \cite{patel2022variationally}), and the VarMiON software \cite{varmion_st_git} is based on the open-source deep learning library PyTorch \cite{pythorch_doc}.
\new{In Table \ref{tab_quantitative}, we present a summary of the performance data for each test case.}
\noindent
\begin{table}
\new{
\begin{tabular}{|p{4.5cm}|c|c|c|}
\hline
 & \textbf{Cavity} & \textbf{Past a cylinder} & \textbf{Contraction} \\
\hline
 size of the dataset & $1000$ & $2002$ & $2000$ \\
\hline
 size of the training set & $800$ & $1601$& $1600$\\
\hline
 size of the validation set & $100$& $200$ & $200$\\
\hline
 size of the test set & $100$ & $200$ & $200$ \\
\hline
sampling strategy for params. & random.uniform & random.uniform& random.uniform \\
\hline
output points in the loss & $225 \times 20$ & $496 \times 20$& $325 \times 20$ \\
\hline
optimizer & torch.optim.Adam &  torch.optim.Adam & torch.optim.Adam\\
\hline
learning rate &$ 10^{-4}$& $ 10^{-4}$&$10^{-4}$ \\
\hline
pdes per batch & $8$ & $8$ & $8$ \\
\hline
\end{tabular}}
\caption{Methods and parameters of the learning procedure for the three different cases: cavity flow, flow past a cylinder and contraction flow.}
\label{tab_details}
\end{table}
\begin{table}
\new{
    \begin{tabular}{|p{3.5cm}|p{2.5cm}|p{2.3cm}|p{2.8cm}|p{2.5cm}|}
        \hline
        & offline data-generation cost & training time & inference time per new parameter instance & speed-up w.r.t. FEM solver \\
        \hline
        Cavity flow & $4800 s$ & $16482.95 s$ & $0.174055 s$ & $0.035$ \\
        \hline
        Flow past a cylinder & $48120 s$ & $94830.72 s$ & $27.277750 s$ &  $1.13$\\
        \hline
        Contraction flow  & $55860 s$ & $57508.50 s$&  $20.078824 s$ &  $0.71$ \\
        \hline
    \end{tabular}
    }
    \caption{Performance data for the three different cases: cavity flow, flow past a cylinder and contraction flow.}
    \label{tab_quantitative}
\end{table}
\subsection{Cavity flow}
The first test is the lid-driven cavity flow, in which we solve the differential problem \eqref{eq:stokes} in a square of width $L$, assuming vanishing initial conditions for the velocity and pressure fields. Moreover, no-slip boundary conditions are imposed on the velocity at the left, bottom and right boundaries (see Figure~\ref{fig:mesh}(a)), while the flow is generated by the horizontal velocity profile
\begin{equation}
\vc u(\vc x, t) = (8 f  (1 + \tanh(8  (t - 0.5))) x_1^2  (1 - x_1)^2, 0)^\top
\end{equation} 
on the top boundary.
Here, the external forcing is represented by the scalar coefficient $f$ in the expression above, while no external force is present in the bulk equation.
The time interval $[0,\tau]$ with $\tau=2$ s is discretized with time step $\dt =10^{-3}$ s.
The training of the network is done with $1000$ instances, each with different values of $\mu\in(1.1,10)$ and $f\in(0.02,0.9)$, \new{that entail $Re = \frac{\rho {{u}_{max}}L }{ \mu} \in (0.002, 0.81)$, with ${u}_{max}$ the maximum speed of the fluid}, for $1000$ cross-validation iterations, named epochs.
The global VarMiON architecture in this case entails the training of about $10^7$ parameters.

\begin{figure}
    \centering
    \includegraphics[ height=0.31\linewidth]{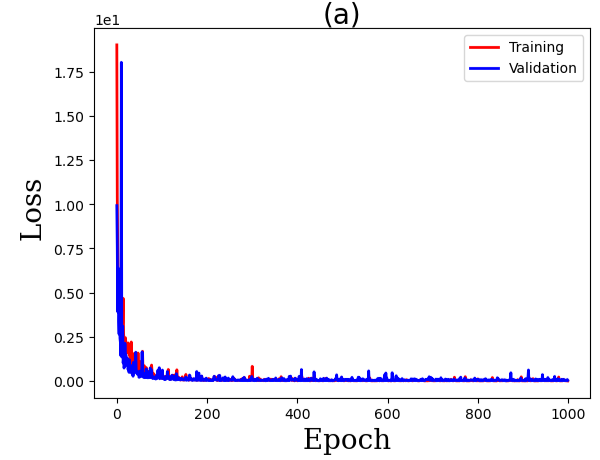}
     \includegraphics[ height=0.31\linewidth]{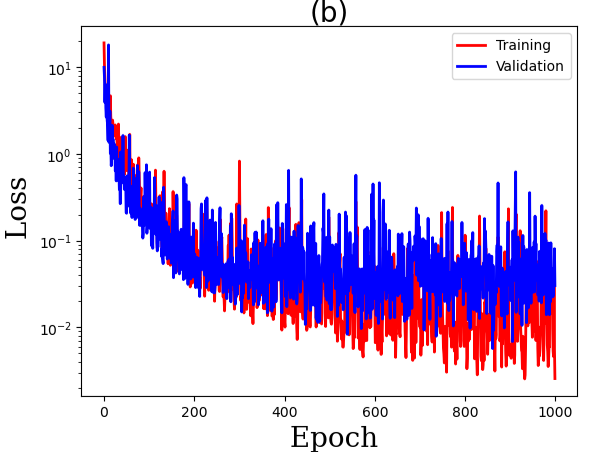}

     \includegraphics[ height=0.33\linewidth]{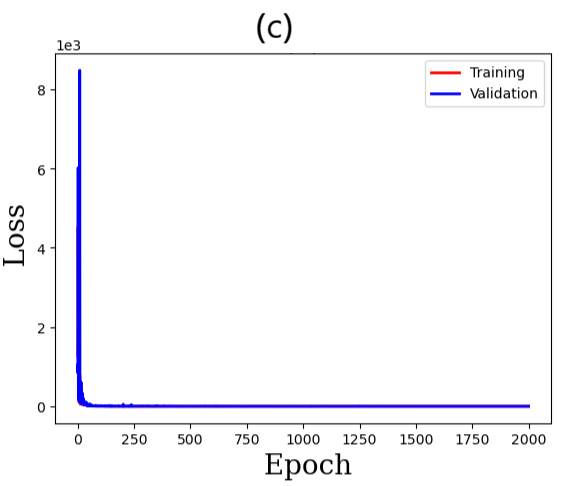}
     \includegraphics[ height=0.34\linewidth]{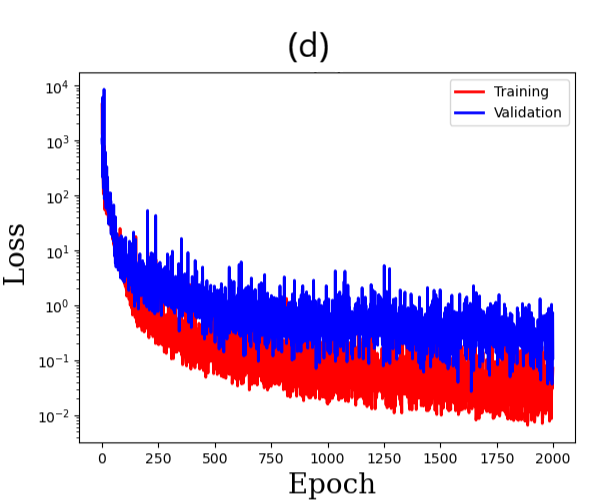}

     \includegraphics[ height=0.33\linewidth]{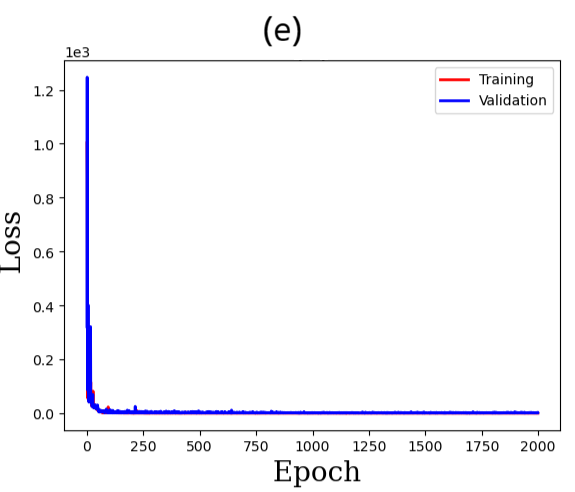}
     \includegraphics[ height=0.33\linewidth]{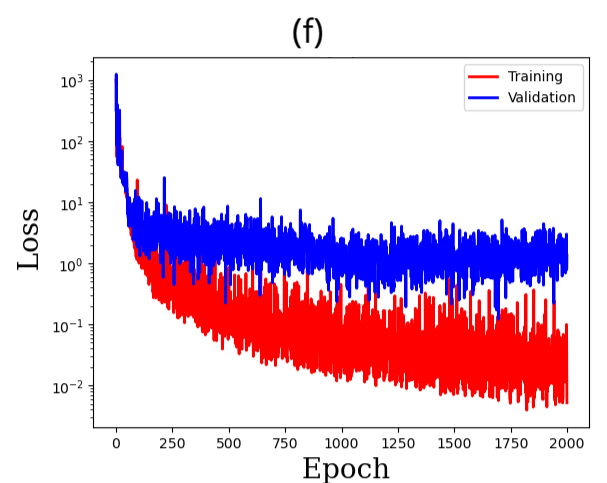}
     
    \caption{The loss function computed, \new{for both velocity and pressure as defined in \eqref{error_l2}}, on the training set (red line) and the validation one (blue line), respectively. Data are presented in linear scale in the left panels and in semi-log scale in the right panels.
    The panels are related to the cavity flow (a)-(b), the flow past a cylinder (c)-(d), and the contraction flow (e)-(f).}
    \label{fig:losesse}
\end{figure}

\begin{figure}
    \includegraphics[height=0.27\linewidth]{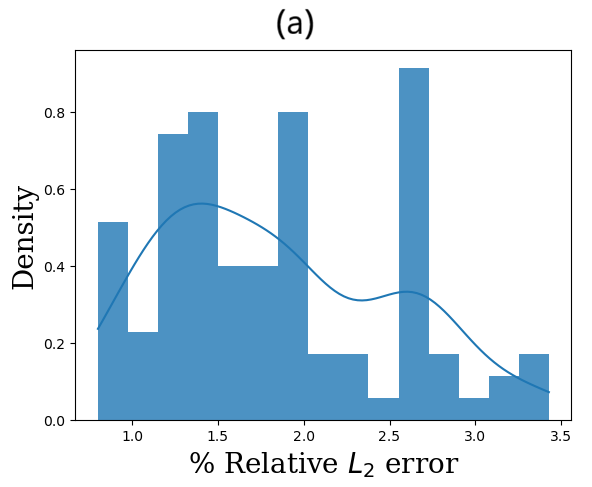}
    \includegraphics[height=0.27\linewidth]{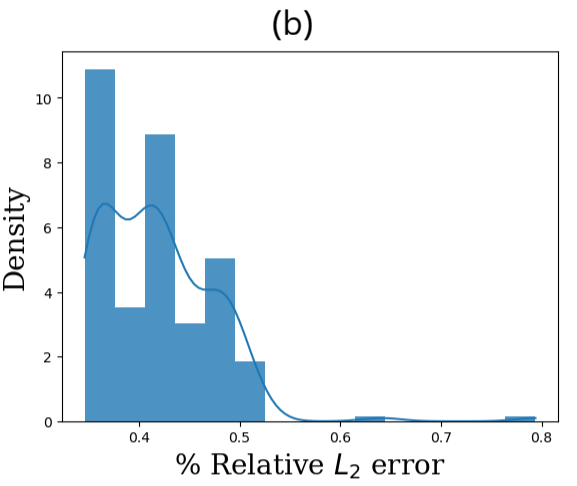}
    \includegraphics[height=0.27\linewidth]{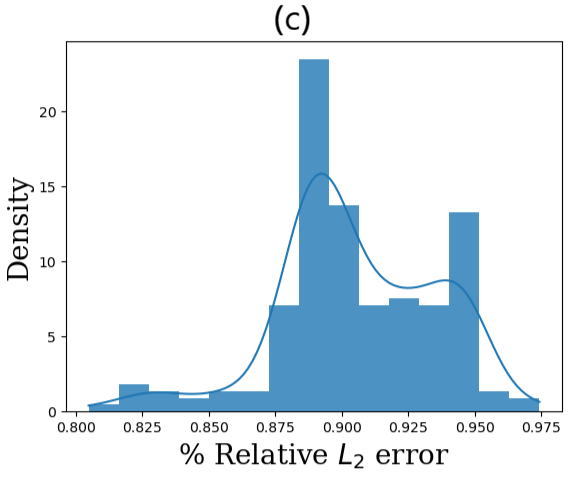}
    
    \caption{The probability density of the relative $L_2$-error computed on the testing set \new{for both velocity and pressure according to  \eqref{error_l2}}. The panels are related to the cavity flow, the flow past a cylinder, and the contraction flow, on the left (a), center (b) and right (c), respectively.}
    \label{fig:probability}
\end{figure}

\begin{figure}
    \includegraphics[height=0.59\linewidth]{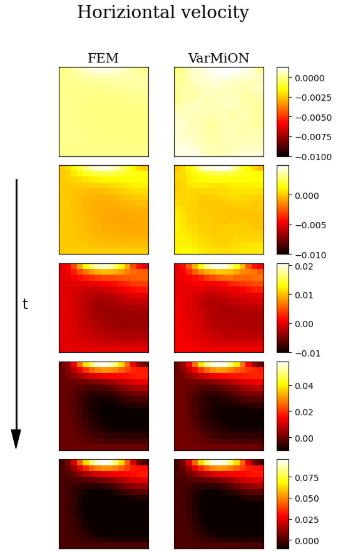}
    \includegraphics[height=0.59\linewidth]{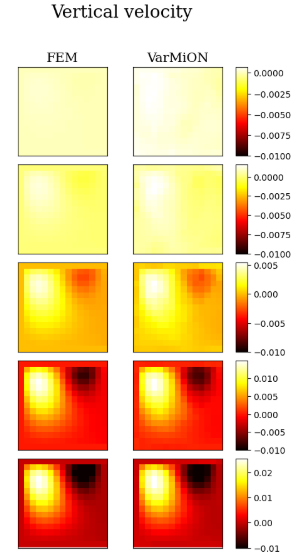}
    \includegraphics[height=0.59\linewidth]{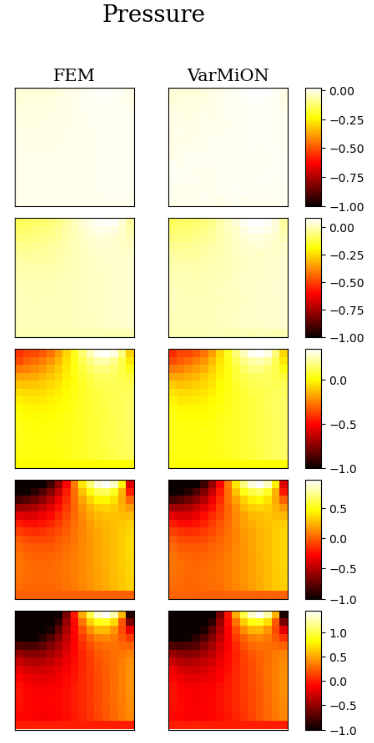}
    \caption{The comparison between the FEM solution and VarMiON prediction \new{of the cavity flow} for the velocity along $x$ and $y$ directions and pressure, on left, center and right, respectively. The plots correspond to $t=0.21, 0.31, 0.42, 0.52, 0.63$ s. \new{The reported test is related to simulations done for $\mu= 2.3456$ and $f = 0.0075$.}}
    \label{fig:comparison_cavity}
\end{figure}

\begin{figure}
    \includegraphics[height=0.26\linewidth]{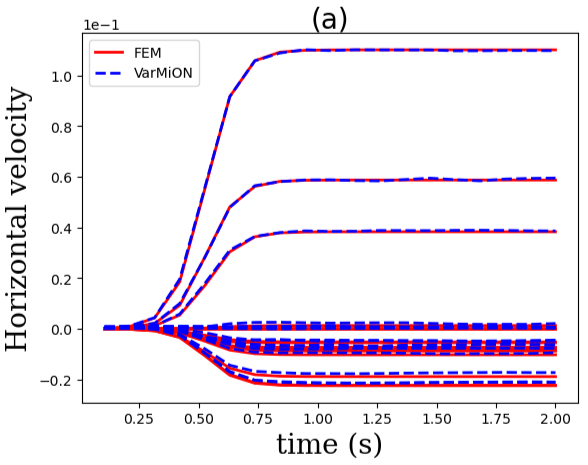}
    \includegraphics[height=0.25\linewidth]{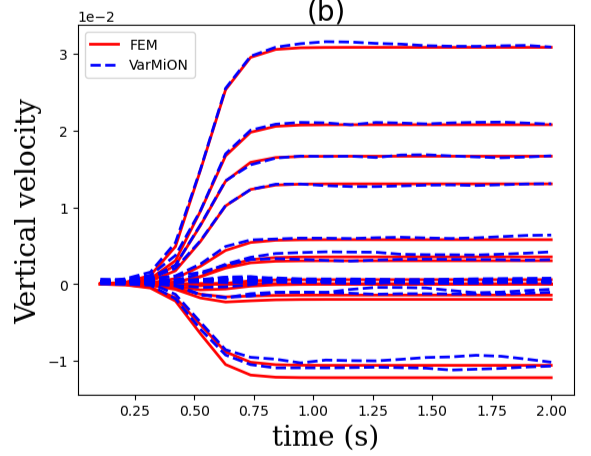}
    \includegraphics[height=0.25\linewidth]{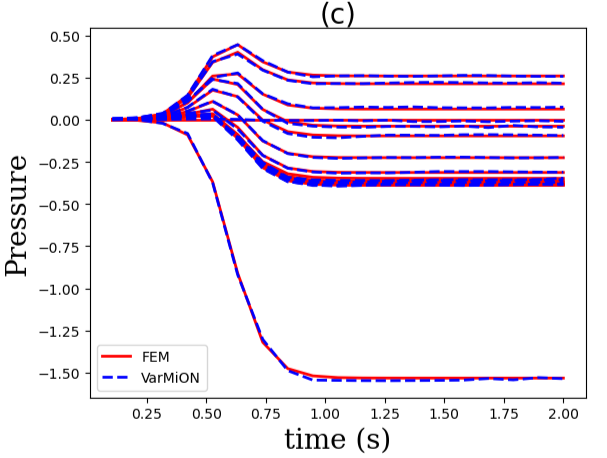}

    \includegraphics[height=0.26\linewidth]{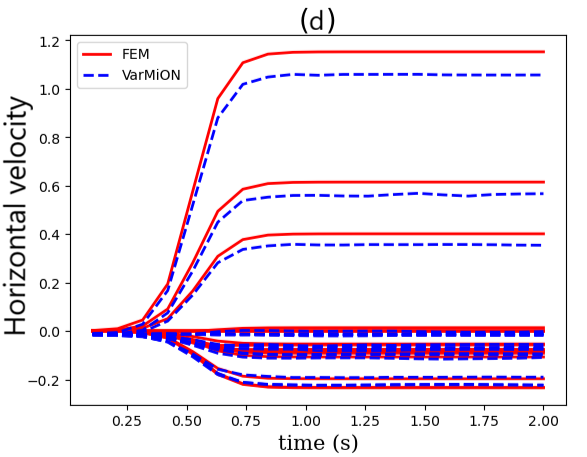}
    \includegraphics[height=0.26\linewidth]{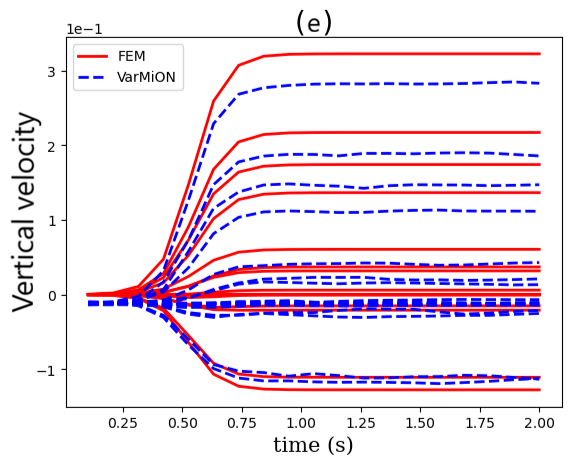}
    \includegraphics[height=0.26\linewidth]{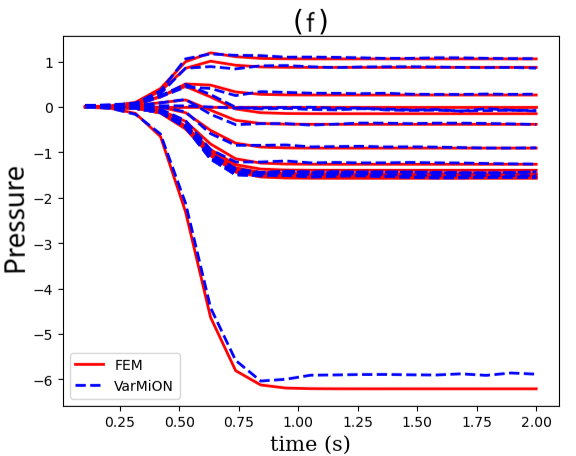}
    \caption{The trend of some chosen points of FEM solution and VarMiON prediction \new{of the cavity flow} for the velocity along $x$ and $y$ directions and the pressure, on left, center, and right, respectively. \new{The tests are related to simulations done with $\mu= 2.3456$ and $f = 0.0075$ (top panels)} \new{and with $\mu= 0.9$ and $f = 1.2$ (bottom panels). The latter values lie slightly outside the training ranges, indeed the predictions are somewhat less accurate, though they remain satisfactory.}}
    \label{fig:trend_cavity}
\end{figure}

After the training, the relative error obtained on the test dataset has a mean of $1.85 \%$ and standard deviation of $0.67 \%$.
\new{Note that the relative error referred to above is the $L_2$ norm calculated simultaneously for both velocity and pressure, as given in formula \eqref{error_l2}.}
\new{Although the datasets differ, it is still meaningful to compare the relative errors of our extended VarMiON with those reported for cavity flow using DeepONet or FNO (and their variants) in \cite{lu2022comprehensive}. VarMiON attains improved accuracy, yielding a mean error of $0.67\%$, which is lower than that of the previous methods ($1.78 \pm 0.22\% - 2.24 \pm 0.21\%$). Another perspective on model performance can be obtained by comparing against the results in \cite{luo2023cfdbench}, where they evaluate an FFN (Mean Square Error $0.78$) and DeepONet (Mean Square Error $4.008$).
In this experiment, the inference time of VarMiON \ref{tab_quantitative} is shorter than that of the  DeepONet--vanilla reported in \cite{lu2022comprehensive}. }
In Figures \ref{fig:losesse}(a)-(b), the losses are shown to reach a steady condition, thereby proving that the number of epochs in the training is sufficient to achieve the desired precision. In Figure \ref{fig:probability}(a) we can see the probability distribution of the relative $L_2$-error, which highlights a sufficiently small spread of the error values.
The comparison between the FEM solution and the VarMiON prediction, presented in Figures \ref{fig:comparison_cavity} and \ref{fig:trend_cavity}, demonstrates an excellent agreement between the numerical and predicted fields. This consistency is evident both in the spatial distribution at different time snapshots (Figure \ref{fig:comparison_cavity}) and in the temporal evolution at selected monitoring points (Figure \ref{fig:trend_cavity}).

\subsection{Flow past a cylinder}
The second test consists of the flow past a cylinder,
in which we solve the differential problem \eqref{eq:stokes} assuming vanishing initial conditions for the velocity and pressure fields. Moreover, no-slip boundary conditions are imposed on the velocity at the top and bottom boundaries and on the cylinder (see Figure~\ref{fig:mesh}(b)). 
The time interval $[0,\tau]$ with $\tau=0.05$ s is discretized with time step $\dt =2.5 \cdot 10^{-5}$ s.
The training of the network is done with $2000$ instances, each with different values of $\mu\in(1.1,10)$. In this case, \new{$\vc f$ is the vector of mean horizontal pressure gradient }$\vc f =(f_x, 0, 0) $ where $f_x\in(80, 720)$. This training is done for $2000$ epochs.
 \new{Here we have $Re =\frac{\rho f_x H^3 }{ 4\mu^2}
  \in (0.0128, 9.520)$, with $H$ the width of the channel.}
The global VarMiON architecture in this case entails the training of about $10^6$ parameters.

\begin{figure}
    \centering
    \includegraphics[width=0.9\linewidth]{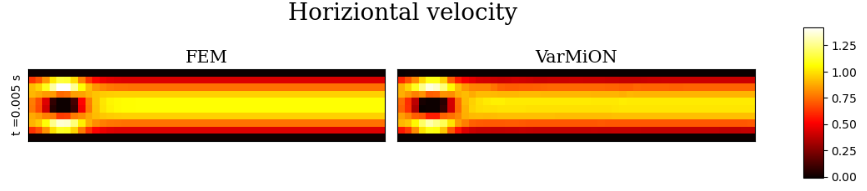}
    \includegraphics[width=0.9\linewidth]{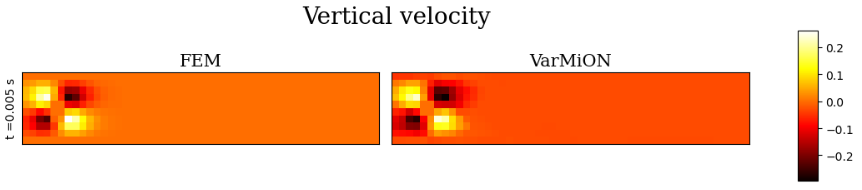}
    \includegraphics[width=0.9\linewidth]{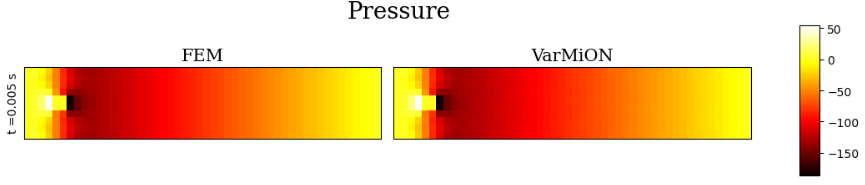}
    \caption{The comparison between the FEM solution and VarMiON prediction at $t=0.005$ s \new{of the flow past a cylinder} for the velocity along $x$ and $y$ directions and the pressure, on top, center and bottom, respectively. The pressure here is the perturbation of the linearly decreasing pressure imposed as a forcing term. \new{The test is related to simulations done with $\mu= 2.712$ and $f_x =304.358$.}}
    \label{fig:comparison_fpc}
\end{figure}

\begin{figure}
    \centering
    \includegraphics[height=0.25\linewidth]{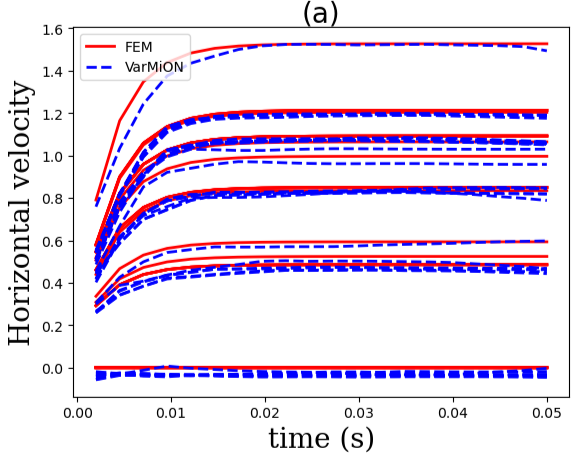}
    \includegraphics[height=0.25\linewidth]{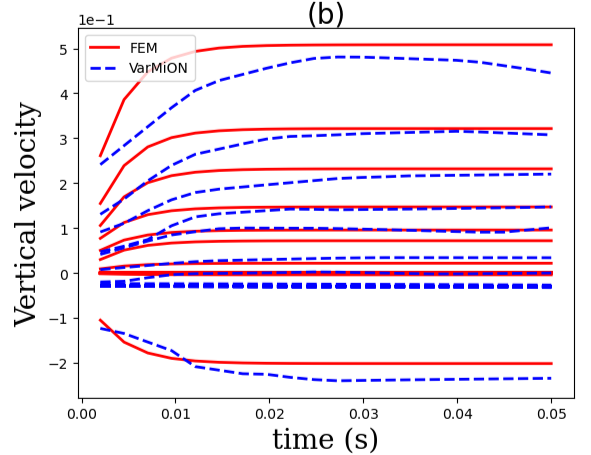}
    \includegraphics[height=0.25\linewidth]{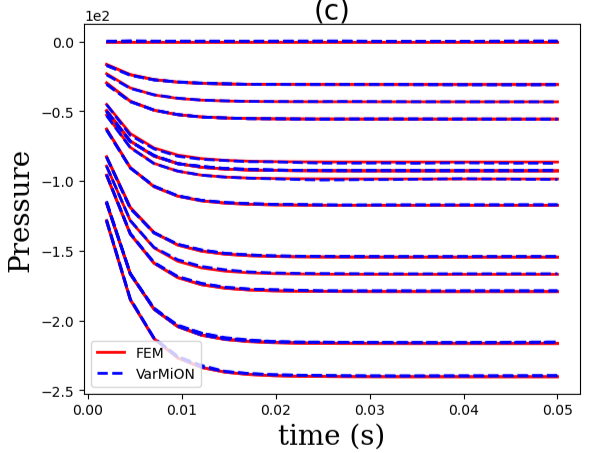}

    \includegraphics[height=0.25\linewidth]{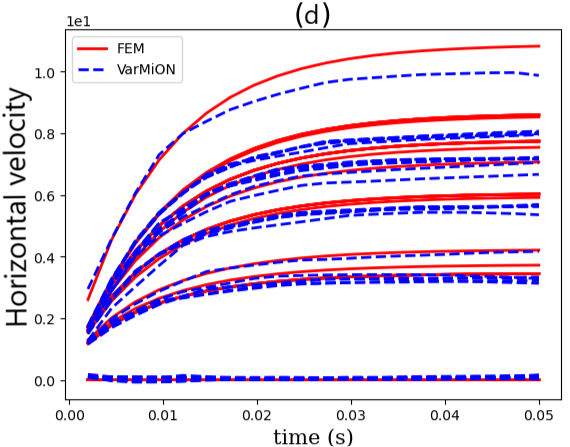}
    \includegraphics[height=0.26\linewidth]{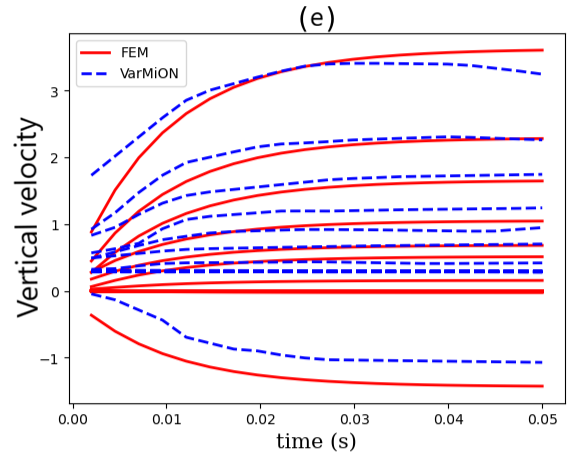}
    \includegraphics[height=0.25\linewidth]{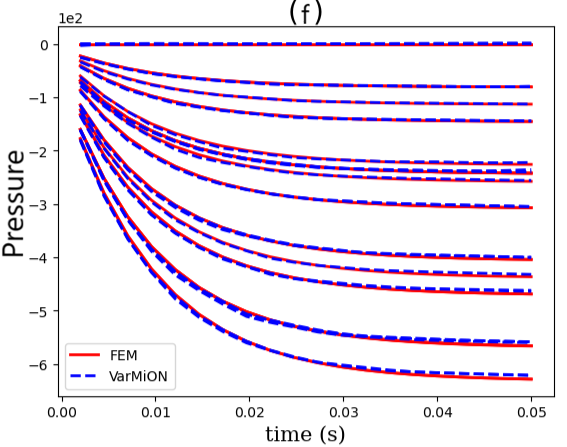}
    \caption{The trend of some chosen points of FEM solution and VarMiON prediction  \new{of the flow past a cylinder} for the velocity along $x$ and $y$ directions and the pressure, on left, center, and right, respectively. The pressure here is the perturbation of the linearly decreasing pressure imposed as a forcing term. \new{The tests are related to simulations done with $\mu= 2.712$ and $f_x =304.358$ (top panels)} \new{and with $\mu= 1$ and $f_x = 800$ (bottom panels). The latter values lie slightly outside the training ranges, indeed the predictions are somewhat less accurate, though they remain satisfactory.}}
    \label{fig:trend_fpc}
\end{figure}

The relative error obtained on the test dataset has a mean of $0.42 \%$ and a standard deviation of $0.06 \%$.
In Figures \ref{fig:losesse}(c)-(d), the losses are shown again to reach a steady condition. 
\new{The relative error obtained with our extended VarMiON, $0.06\%$, is smaller than that reported in the numerical example of \cite{shinde2016galerkin} using a standard ROM-POD approach, which is $0.25 \%$. The performance can also be assessed by comparing with the results in \cite{luo2023cfdbench}, where they analyze FFN (Mean Square Error $0.047$) and DeepONet (Mean Square Error $3.032$).}
In Figure \ref{fig:probability}(b) we can see the probability distribution of the relative $L_2$-error, which highlights a sufficiently small spread of the error values.
As for the cavity flow example, it is evident the good agreement between the FEM solution and the VarMiON in both the spatial distribution (Figure \ref{fig:comparison_fpc}) and the temporal evolution (Figure \ref{fig:trend_fpc}).

\subsection{Contraction flow}
The third test consists of the contraction flow in which we solve the differential problem \eqref{eq:stokes} assuming vanishing initial conditions for the velocity and pressure fields. As usual, we simulate only half of the domain, imposing no-slip boundary conditions on the top boundary profile, while the vertical velocity is set to zero on the bottom boundary, as required by the symmetry of the entire flow (see Figure~\ref{fig:mesh}(c)).
The time interval $[0,\tau]$ with $\tau=0.05$ s is discretized with time step $\dt =2.5 \cdot 10^{-5}$ s.
The training of the network is done with $2000$ instances, each with different values of $\mu\in(1.1,10)$. In this case, \new{$\vc f$ is the vector of mean horizontal pressure gradient} $\vc f =(f_x, 0, 0) $ where $f_x\in(80, 720)$. This training is done for $2000$ epochs.
 \new{Here we have $Re = \frac{\rho f_x H^3 }{ 4\mu^2}\in (0.0128, 9.520)$, with $H$ the width of the inlet.}
The global VarMiON architecture in this case entails the training of about $10^6$ parameters.

The relative error obtained on the test dataset has a mean of $0.91 \%$ and a standard deviation of $0.03 \%$.
\new{Although the datasets are not identical, it is still informative to compare the relative errors of our extended VarMiON with the results reported in \cite{demo2023deeponet} for POD, DeepONet, and their MFDeepONet method. VarMiON attains a comparable accuracy, with a mean error of $0.91\%$, while their MF approach reaches a best-case error of $0.9\%$.}

\begin{figure}
    \centering
    \includegraphics[width=0.9\linewidth]{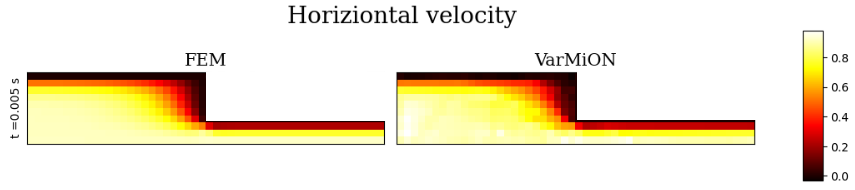}
    \includegraphics[width=0.9\linewidth]{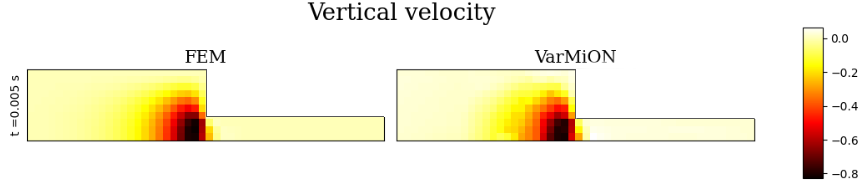}
    \includegraphics[width=0.9\linewidth]{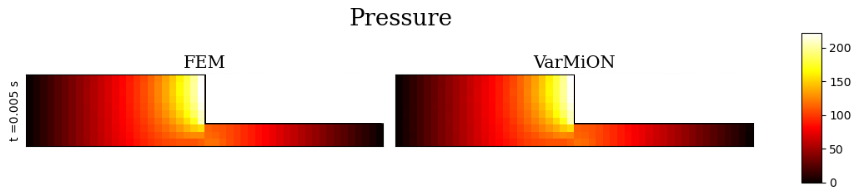}
    \caption{The comparison between the FEM solution and VarMiON prediction at $t=0.005$ s \new{of the contraction flow} for the velocity along $x$ and $y$ directions and the pressure, on top, center and bottom, respectively. The pressure here is the perturbation of the linearly decreasing pressure imposed as a forcing term. \new{The test is related to simulations done with $\mu= 1.6434$ and $f_x = 445.305$.}}
    \label{fig:comparison_cf}
\end{figure}

\begin{figure}
    \includegraphics[height=0.256\linewidth]{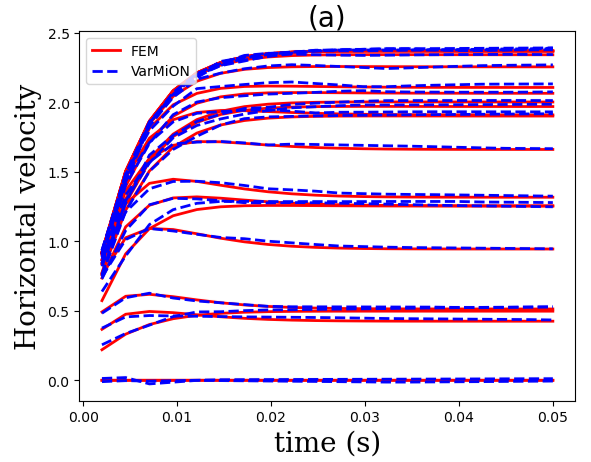}
    \includegraphics[height=0.256\linewidth]{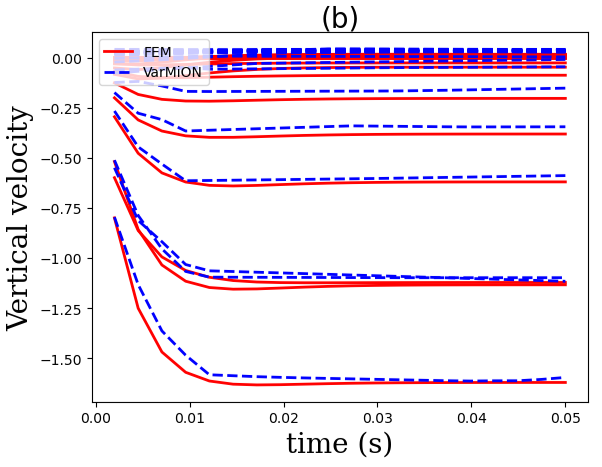}
    \includegraphics[height=0.26\linewidth]{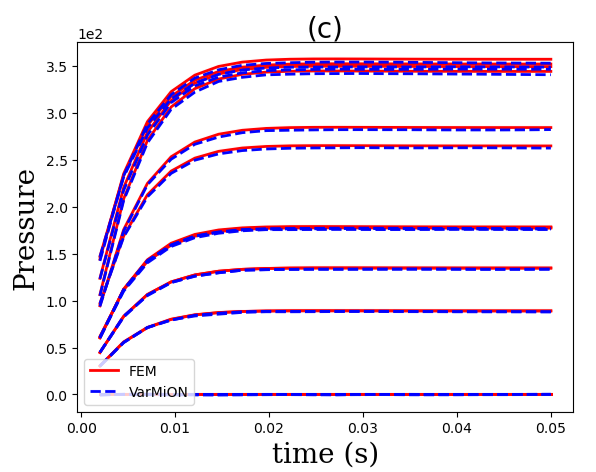}

    \includegraphics[height=0.26\linewidth]{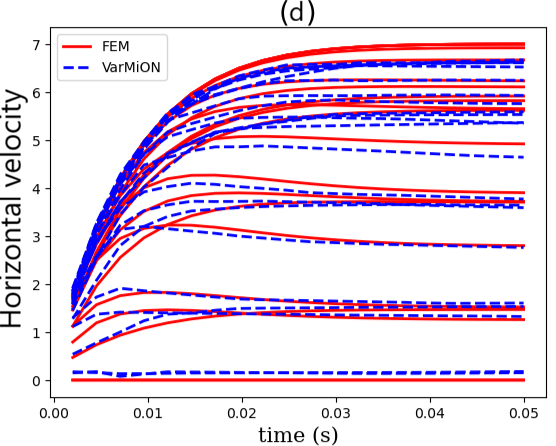}
    \includegraphics[height=0.26\linewidth]{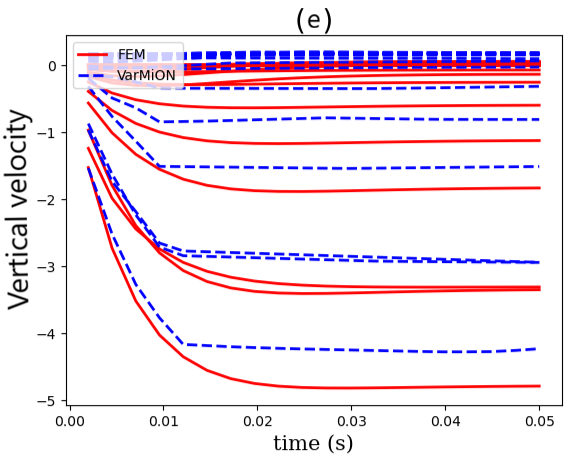}
    \includegraphics[height=0.26\linewidth]{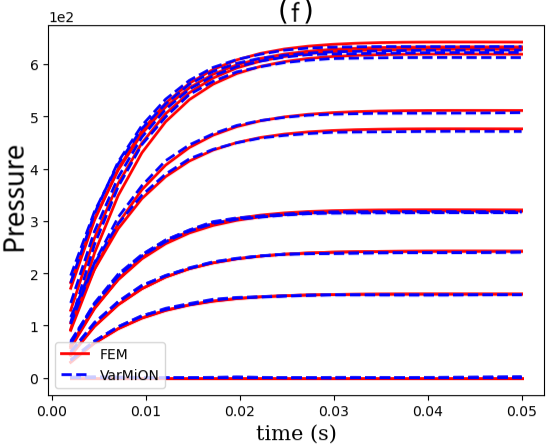}
    \caption{The trend of some chosen points of the FEM solution and VarMiON prediction \new{of the contraction flow}  for the velocity along $x$ and $y$ directions and the pressure, on the left, center, and right, respectively. The pressure here is the perturbation of the linearly decreasing pressure imposed as a forcing term. \new{The tests are related to simulations done with $\mu= 1.6434$ and $f_x = 445.305$ (top panels)} \new{and with $\mu =1$ and $f_x = 800$ (bottom panels). The latter values lie slightly outside the training ranges, indeed the predictions are somewhat less accurate, though they remain satisfactory.}}
    \label{fig:trend_cf}
\end{figure}

In Figures \ref{fig:losesse}(e)-(f), the losses are shown to reach a steady condition. In Figure \ref{fig:probability}(c) we can see the probability distribution of the relative $L_2$-error.
Once again, there is a strong concordance between the FEM solution and the VarMiON results, both in the spatial distributions (Figure \ref{fig:comparison_cf}) and in the temporal evolution ({Figure} \ref{fig:trend_cf}).

\begin{figure}
    \centering
    \includegraphics[width=0.9\linewidth]{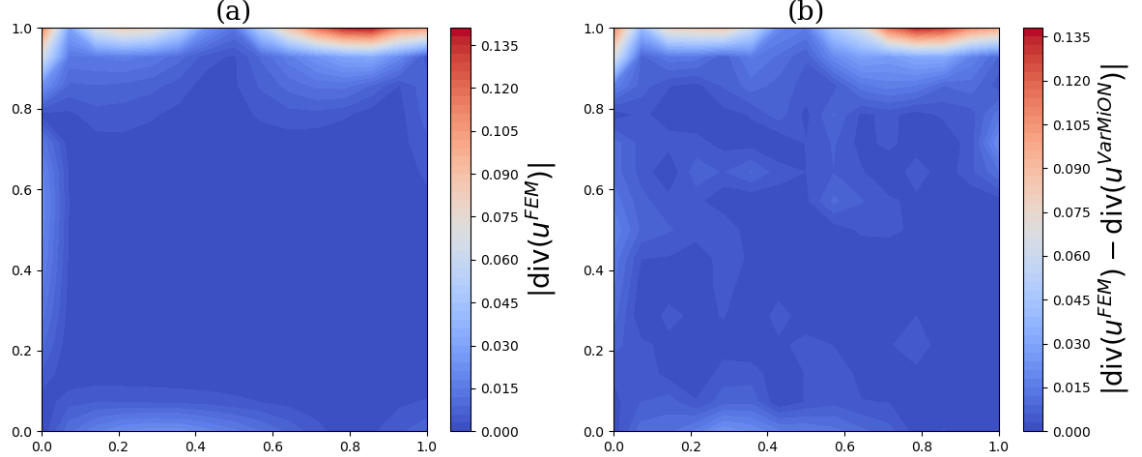}
    
    \includegraphics[width=0.9\linewidth]{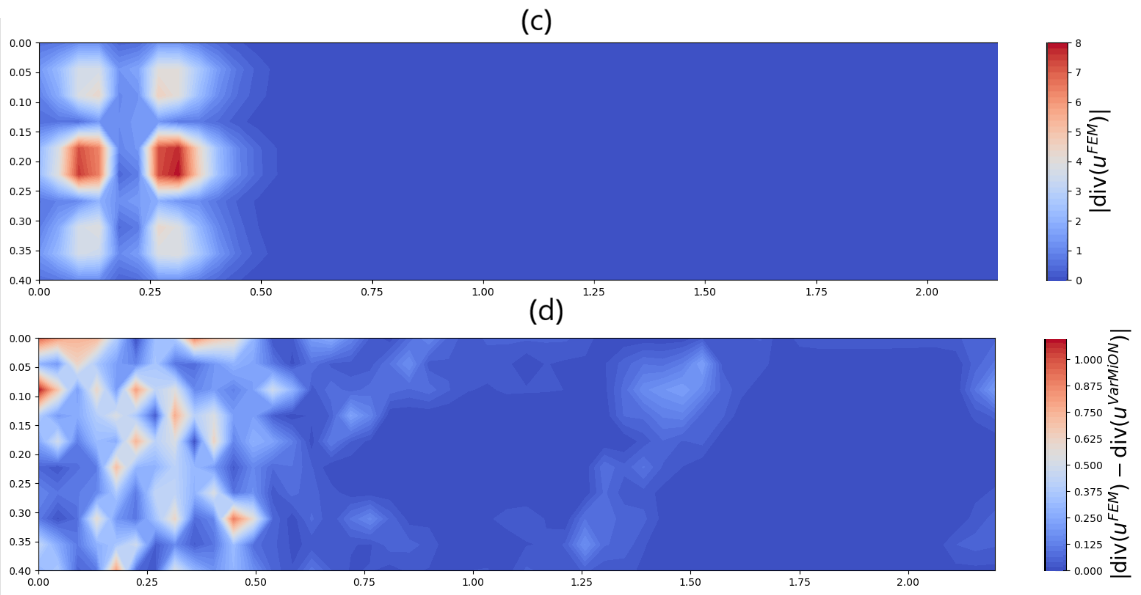}
    
    \includegraphics[width=0.9\linewidth]{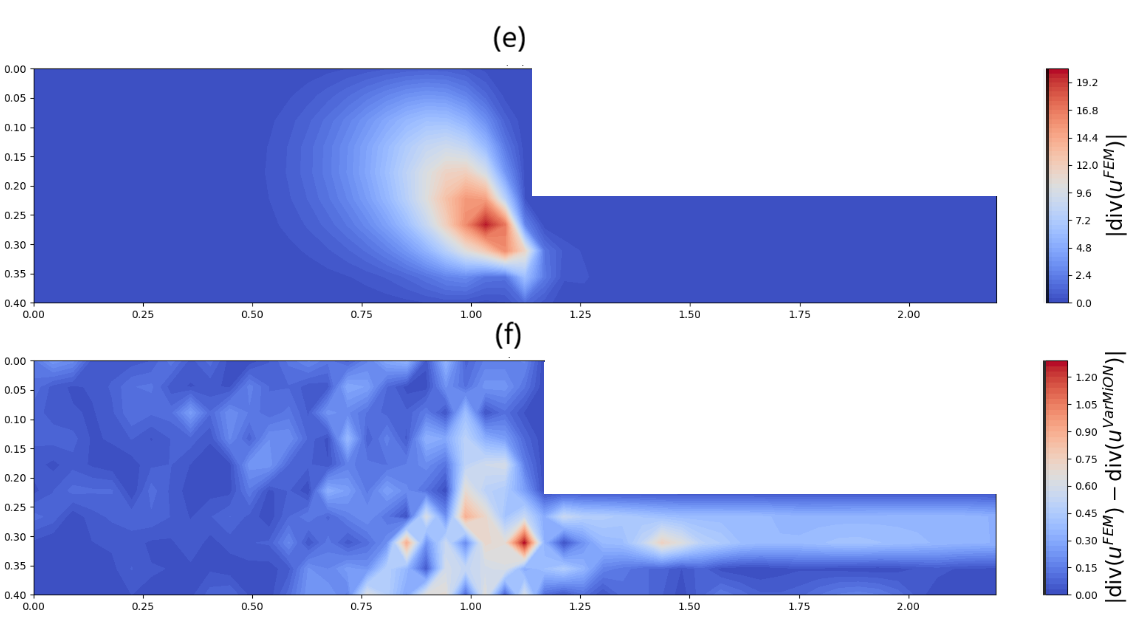}
    \caption{\new{ The divergence of the velocity field $\vc u^{FEM}$ and the absolute error between the divergence of $\vc u^{FEM}$ and $\vc u^{VarMiON}$ for the cavity flow (a)-(b), the flow past a cylinder (c)-(d), and the contraction flow (e)-(f). The VarMiON respects the incompressibility constraint with an accuracy comparable with the one present in the training set.}}
    \label{fig:div}
\end{figure}

\new{As a final remark, relevant to all three of the geometries that we considered, we discuss to what extent the incompressibility constraint is respected.
To do this, we evaluate the divergence of the predicted velocity field, that should be zero in the exact solution, and compare it against the finite-element reference. In our approach, the divergence-free condition is enforced in a weak sense. We observe that the VarMiON prediction captures this constraint as accurately as the reference FEM method (Figure \ref{fig:div}). 
}

\section{Conclusions}\label{sec:conclusions}

We have presented a study devoted to forecasting the velocity and pressure fields governed by the time-dependent Stokes equations. 
We employ an operator-learning neural network that explicitly incorporates the variational formulation of the problem to approximate its solution. This approach can be interpreted as a hybrid methodology between the classical DeepONet and Physics-Informed Neural Networks (PINNs), as it embeds the underlying physical principles directly into the network architecture, thereby enhancing the interpretability of the model.

Within this framework, our contribution consists of extending the VarMiON, originally proposed for elliptic PDEs, to the prediction of the time-dependent velocity and pressure fields arising in the unsteady Stokes problem.
Our numerical experiments demonstrate the effectiveness of this machine learning strategy in a fluid mechanics context.  

The numerical results indicate that, for applications in fluid dynamics, VarMiON represents a competitive alternative to more traditional neural network approaches. In particular, the $L_2$-error between the reference numerical solution and the  prediction is sufficiently small.
Moreover, the temporal and spatial evolution of the velocity and pressure fields is accurately captured by the VarMiON-based predictions.

The present application of the VarMiON approach is intended as an initial step toward the application of this method for fluid dynamics problems. 
In future works, we plan to investigate its performance at higher Reynolds numbers, where conventional numerical simulations become computationally demanding and a VarMiON-based strategy may offer significant advantages.
A further interesting research direction would be the application of VarMiON to more complex flows of non-Newtonian fluids.


\section*{Code availability}

The code used in the study is publicly available from the GitHub repository 

\noindent \url{https://github.com/NLADlab/VarMiON/Stokes}.


\section*{Acknowledgements}
The work of L.R. is supported by the Italian Ministry of Research, under the complementary actions to the
NRRP “D34Health - Digital Driven Diagnostics, prognostics and therapeutics for sustainable Health
care” Grant \#PNC0000001, CUP B53C22006100001 and is partially supported by the ``Gruppo Nazionale per il Calcolo Scientifico'' (GNCS - INdAM). L.R. acknowledges the support of the Department of Mathematics of the University of Padua through a research fellowship
during which part of this work was developed.

The work of G.G.G. is partially supported by the ``Gruppo Nazionale per la Fisica Matematica'' (GNFM - INdAM).

The authors would like to thank Fabio Marcuzzi for several useful discussions.

\bibliographystyle{abbrvurl}
\bibliography{Bibliography}

\end{document}